\documentclass[12pt,draft]{amsart}
\usepackage{amsmath,amsthm,latexsym,amscd,amsbsy,amssymb,fleqn,leqno}
\chardef\bslash=`\\ 

\makeatletter
\def\verbatim{\interlinepenalty\@M \@verbatim
  \leftskip\@totalleftmargin\advance\leftskip2pc
  \frenchspacing\@vobeyspaces \@xverbatim}
\makeatother
\makeatletter
  \def\dgt@k{\dg@DX=-3 \dg@DY=2 \dg@SIZE=3}
\makeatother

\makeatletter
  \def\dgt@kk{\dg@DX=3 \dg@DY=-1 \dg@SIZE=3}%
\makeatother

\theoremstyle{plain}
\newtheorem{thm}{Theorem}[section]
\newtheorem{cor}[thm]{Corollary}
\newtheorem{lem}[thm]{Lemma}
\newtheorem{pro}[thm]{Proposition}

\theoremstyle{definition}

\numberwithin{equation}{section}

\newcounter{rmnum}

\def\symbolnote#1#2{\let\thefootn=\thefootnote%
\renewcommand{\thefootnote}{\fnsymbol{footnote}}%
\footnotemark[#1]%
\footnotetext[#1]{#2}%
\let\thefootnote=\thefootn
}

\newfont{\bbb}{msbm10 scaled \magstep1}
\newfont{\bbc}{msbm8 scaled \magstep0}

\newcommand{\R}{\mbox{\bbb R}}

\newcommand{\N}{\mbox{\bbb N}}
\newcommand{\Q}{\mbox{\bbb Q}}

\newcommand{\mesh}{\hbox{\rm mesh}}

\begin{document}


\title{ Roberts' type embeddings and conversion of
the transversal Tverberg's theorem}
\author{S.A. Bogatyi}
\address{Faculty of Mechanics and Mathematics,
Moscow State University, Vorob'evy gory, Moscow, 119899 Russia}
\email{bogatyi@mech.math.msu.su}
\thanks{The paper was started during the first author's visit
to Nipissing University in May 2004.
He thanks NSERC for supporting his visit and COMA Department of NU
for hospitality. The first author also acknowledges  a financial 
support from RFFI (grant 03-01-00706).}
\author{V.  Valov}
\address{Department of Computer Science and Mathematics, Nipissing University,
100 College Drive, P.O. Box 5002, North Bay, ON, P1B 8L7, Canada}
\email{veskov@nipissingu.ca}
\thanks{The second author was partially supported by NSERC Grant 261914-03}

\keywords{embedding, covering dimension, finite-to-one maps, $k$-regular maps,
Tverberg's transversal theorem}
\subjclass{Primary: 54F45; Secondary: 55M10, 54C65.}


\begin{abstract}{Here are two of our main results:\\
{\bf Theorem 1.}
Let $X$ be a normal space with $\dim X=n$ and $m\geq n+1$.
Then the space $C^*(X,\R^m)$ of all bounded maps from $X$ into $\R^m$
equipped with the uniform convergence topology contains a dense
$G_{\delta}$-subset consisting of maps $g$ such that
$\overline{g(X)}\cap\Pi^d$ is at most $(n+d-m)$-dimensional for every
$d$-dimensional plane $\Pi^d$ in $\R^m$, where $m-n\leq d\leq m$.

\smallskip\noindent
{\bf Theorem 2.}
Let $X$ be a metrizable compactum with $\dim X\leq n$ and $m\geq n+1$.
Then, $C(X,\R^m)$ contains a dense $G_{\delta}$-subset of maps $g$ such that
for any integers $t,d,T$ with $0\leq t\leq d\leq m-n-1$ and $d\leq T\leq m$
and any $d$-plane $\Pi^d\subset\R^m$ parallel to some coordinate planes
$\Pi^t\subset\Pi^T$ in $\R^m$, the inverse image $g^{-1}(\Pi^d)$ has at most
$q$ points, where $\displaystyle q=d+1-t+\frac{n+(n+T-m)(d-t)}{m-n-d}$
if $n\geq (m-n-T)(d-t)$ and
$\displaystyle q=1+\frac{n}{m-n-T}$ otherwise.

\smallskip\noindent
In case $m=2n+1$, the combination of Theorem 1 and
the N\"{o}beling--Pontryagin embedding theorem provides a generalization of
a theorem due to Roberts \cite{r}.
Theorem 2 extends the following results:
the N\"{o}beling--Pontryagin embedding theorem ($d=0$, $m=T\geq 2n+1$);
Hurewicz's theorem \cite{h} about mappings into an Euclidean space with
preimages of small cardinality ($d=0$, $n+1\leq m=T\leq 2n$);
Boltyanski's theorem \cite[Theorem 1]{b} about $k$-regular maps
($d=k-1$, $t=0$, $T=m\geq nk+n+k$) and
Goodsell's theorem \cite{g1} about existence of special embeddings
($t=0$, $T=m$).
An infinite-dimensional analogue of Theorem 2 is also established.
Our results are based on Theorem 1.1 below which is considered as
a converse assertion of the transversal Tverberg's theorem and implies
the Berkowitz-Roy theorem \cite{br}, \cite{g1}.}
\end{abstract}

\maketitle
\markboth{S.~Bogatyi and V.~Valov}{Roberts' type embeddings}


\section{Introduction}

\bigskip
All maps considered in our paper are continuous and all spaces are assumed
to be at least completely regular.
Everywhere below $\Pi^k\subset\R^m$ denotes a $k$-dimensional,
not necessarily coordinate, plane (simply, a $k$-plane) in $\R^m$.
Moreover, if $\Pi^t$ and $\Pi^T$ are two coordinate planes with $t\leq T$,
we write $\Pi^t\subset\Pi^T$ if $\Pi^t$ is a linear subspace of $\Pi^T$.
If not explicitely stated otherwise, all function spaces in the paper are
equipped with the source limitation topology.

Our goal is to prove Theorem 1.1 below and provide some applications of
this theorem.

\begin{thm}
Let $\displaystyle A_{i,j}$, $i=1,2,...,q$ and
$\displaystyle j=1,2,...,n_{i}+1$, be points in $\R^m$ such that
the set of their coordinates is algebraically independent.
Suppose $0\leq t\leq d\leq T\leq m$ and $\Pi^d$ is a $d$-plane in $\R^m$
parallel to some coordinate planes $\Pi^t\subset\Pi^T\subset\R^m$.
If either $d-t+1\leq q$ and
$\displaystyle n_1+n_2+...+n_q+1\leq (m-d)(q-1)-(T-d)(d-t)$
or $q\leq d-t+1$ and $\displaystyle n_1+n_2+...+n_q+1\leq (m-T)(q-1)$,
then there exists $i\in\{1,2,...,q\}$ such that $\Pi^d$ doesn't meet
the linear hull $\Pi(M_i)$ of the set
$\displaystyle M_i=\{A_{i,1},...,A_{i,n_{i}+1}\}$.
\end{thm}

The first part of Theorem 1.1, when $t=0$, $T=m$ and $d-t+1\leq q$, was stated
as a conjecture in \cite[Conjecture 2]{b2} and \cite[Conjecture 4.2]{b4}.

Recall that a real number $v$ is called algebraically dependent
on the real numbers $u_1,...,u_k$ if $v$ satisfies the equation
$p_0(u)+p_1(u)v+...+p_n(u)v^n=0$, where $p_0(u),...,p_n(u)$ are polynomials
in $u_1,...,u_k$ with rational coefficients, not all of them 0.
A finite set of real numbers is algebraically independent if none of them
depends algebraically on the others.

\begin{cor}
Let $K$ be a finite simplicial complex, $\theta\colon K\to\R^m$
a semi-linear map and $\varepsilon >0$.
Then there is a semi-linear map $g\colon K\to\R^m$ such that
$d(g(v),\theta(v))<\varepsilon $ for each vertex $v$ of $K$, and for
any integers $n,d,t,T$ with $0\leq t\leq d\leq m-n-1$ and $d\leq T\leq m$,
and any $d$-plane $\Pi^d\subset\R^m$ parallel to some coordinate planes
$\Pi^t\subset\Pi^T\subset\R^m$, the number $q$ of pairwise disjoint simplexes
of $K$ of dimension $\leq n$ whose images under $g$ intersect $\Pi^d$ satisfies
the inequalities:
$\displaystyle q\leq d+1-t+\frac{n+(n+T-m)(d-t)}{m-n-d}$ if $n\geq (m-n-T)(d-t)$
and $\displaystyle q\leq 1+\frac{n}{m-n-T}$ if $n\leq (m-n-T)(d-t)$.
\end{cor}

The idea to use algebraically independent sets for proving general position
theorems, like Corollary 1.2, was originated by Roberts in \cite{r}.
This idea was also applied by Berkowitz and Roy in \cite{br} where they
stated a version of Corollary 1.2 with $t=0$ and $T=m$.
A proof of the Berkowitz-Roy theorem was provided by Goodsell
in \cite[Theorem A.1]{g2}
(see also \cite{g1} for another application of the Berkowitz-Roy theorem).

Here are some applications of Corollary 1.2.

\begin{thm}
Let $f\colon X\to Y$ be a perfect map between paracompact spaces such that
$\dim f\leq n$ and $\dim Y=0$.
Then, for every $m\geq n+1$, $C^*(X,\R^m)$ contains a dense $G_{\delta}$-subset
$\mathcal H$ of maps $g$ such that $\dim g(f^{-1}(y))\cap\Pi^d\leq n+d-m$ for
every $y\in Y$ and every $d$-plane $\Pi^d\subset\R^m$ with $m-n\leq d\leq m$.
\end{thm}

\begin{cor}
Let $X$ be a normal space with $\dim X=n$ and $m\geq n+1$.
Then $C^*(X,\R^m)$ equipped with the uniform convergence topology contains
a dense $G_{\delta}$-subset $\mathcal H$ consisting of maps $g$ such that
$\overline{g(X)}\cap\Pi^d$ is at most $(n+d-m)$-dimensional for every
$d$-plane $\Pi^d\subset\R^m$ with  $m-n\leq d\leq m$.
\end{cor}

Let show how Corollary 1.4 follows from Theorem 1.3.
Considering the \v{C}ech-Stone compactification $\beta X$ of $X$ and the
function space $C(\beta X,\R^m)$ instead, respectively, of $X$ and
$C^*(X,\R^m)$, we can suppose that $X$ is compact.
Then take $f$ to be a constant map on $X$ and apply Theorem 1.3.

Roberts proved in \cite[Theorem 1.2]{r} that if $X$ is a compact metrizable
space of dimension $\leq n$ and $n+1\leq d\leq 2n+1$,
then $C(X,\R^{2n+1})$ with the uniform convergence topology contains a dense
$G_{\delta}$-subset consisting of maps $g$ such that
$\dim g(X)\cap\Pi^d\leq d-n-1$ for every $d$-plane $\Pi^d\subset\R^{2n+1}$.
Using this result and the Hurewicz theorem about metrizable compactifications
preserving dimension, Roberts derived the existence of such embeddings
for separable metrizable spaces of dimension $\leq n$.
Obviously, this results of Roberts follow from the combination of Corollary 1.4
and the N\"{o}beling--Pontryagin embedding theorem.

When $d\leq m-n-1$, Theorem 1.3 doesn't work.
But next theorem shows that, in this case, we can prove something more,
we can find a residual subset of $C^*(X,\R^m)$ consisting of maps $g$
such that $g(f^{-1}(y))\cap\Pi^d$ is finite for any $y\in Y$ and
any $d$-plane in $\R^m$.

\begin{thm}
Suppose $f\colon X\to Y$ is a perfect map between metrizable spaces
with $\dim f\leq n$ and $\dim Y\leq 0$.
Then, $C^*(X,\R^m)$ contains a dense $G_{\delta}$-subset $\mathcal{K}$
of maps $g$ such that, for any integers $d,t,T$ with
$0\leq t\leq d\leq T\leq m$ and $d\leq m-n-1$ and any $d$-plane
$\Pi^d\subset\R^m$ parallel to some coordinate planes
$\Pi^t\subset\Pi^T\subset\R^m$, each set
$f^{-1}(y)\cap g^{-1}(\Pi^d)$, $y\in Y$, has at most $q$ points, where
$\displaystyle q=d+1-t+\frac{n+(n+T-m)(d-t)}{m-n-d}$ if $n\geq (m-n-T)(d-t)$
and $\displaystyle q=1+\frac{n}{m-n-T}$ otherwise.
\end{thm}

Taking $Y$ in Theorem 1.5 to be a point, 
$m=n+2$, $d=1$, $t=0$ and $T=r$ we obtain next corollary.

\begin{cor}
Let $X$ be a metrizable compactum with $\dim X\leq n$.
Then, $C(X,\R^{n+2})$ contains a dense $G_{\delta}$-subset of maps $g$
such that for any natural integer $r$ with $r\leq n+2$ and any
line $\Pi^1\subset\R^{n+2}$ parallel to some coordinate plane $\Pi^r$ in $\R^{n+2}$,
the inverse image $g^{-1}(\Pi^1)$ has at most
$n+r$ points.
\end{cor}

Theorem 1.5 together with \cite[Theorem 5.1]{w} implies the following interesting fact:
If $X$ is an $n$-dimensional metrizable compactum with $n\geq 3$ and $m\geq 2n+1$, then almost every embedding of $X$ in $\R^m$ is tame.

Theorem 1.7 below is an infinite-dimensional version of Theorem 1.5.
For compact $X$ and $Y$ --- a point, Theorem 1.7 was established by
Boltyanski \cite{b} under the additional restriction that $r=0$.

\begin{thm}
Let $f\colon X\to Y$ be a perfect map between metrizable spaces with
$Y$ being a $C$-space.
For any integers $d,r$ let $\mathcal{P}(d,r)$ denote the family of all
$d$-planes $\Pi^d\subset l_2$ parallel to some coordinate plane
$\Pi^r\subset l_2$.
Then, $C^*(X,l_2)$ contains a dense $G_{\delta}$-subset of maps $g$ such that,
for every $y\in Y$, and every $\Pi^d\in \mathcal{P}(d,r)$,
$f^{-1}(y)\cap g^{-1}(\Pi^d)$ has at most
$\displaystyle d+1-r$ points if $r\leq d$ and at most one point if $r\geq d$.
\end{thm}

The paper is organized as follows.
The proofs of Theorem 1.1 and Corollary 1.2 are given in Section 2.
Section 3 is devoted to the proof of Theorem 1.3.
Theorems 1.5 and 1.7 are established in Section 4.
The final Section 5 contains more applications of Theorem 1.1 and Corollary 1.2.
Some conjectures are also included in the final section.

A few words about the source limitation topology. 
For any spaces $M$ and $K$ by $C(K,M)$ we denote the set of all continuous maps
from $K$ into $M$.
If $(M,d)$ is a metric space and $K$ is any space, then
the source limitation topology on $C(K,M)$ is defined in the following way:
a subset $U\subset C(K,M)$ is open in $C(K,M)$ with respect to
the source limitation topology provided for every $g\in U$ there exists
a continuous function $\alpha\colon K\to (0,\infty)$ such that
$\overline{B}(g,\alpha)\subset U$.
Here, $\overline{B}(g,\alpha)$ denotes the set
$\{h\in C(K,M):d(g(x),h(x))\leq\alpha (x)\hbox{}~~\mbox{for each
$x\in K$}\}$.
It is well known (see, for example \cite{jm:75}) that if
$(M,d)$ is a complete metric space, $C(K,M)$ with this topology
has Baire property.
This implies that $C^*(K,H)$ with the source limitation topology
also possesses the Baire property for any Banach space $H$.

In conclusion of the introduction, we wish to thank Prof. T. Goodsell
who provided us with his proof \cite[Appendix]{g2}of the Berkowitz-Roy theorem.



\section{Proof of Theorem 1.1 and Corollary 1.2}

{\em Proof of Theorem $1.1$.}
Suppose $\Pi^d$ meets the linear hull $\Pi(M_i)$ of each set $M_i$ and let
$Y_i\in\Pi^d\cap\Pi(M_i)$, $i=1,...,q$.
It suffices to show that, under this assumption, we have either
$\displaystyle n_1+n_2+...+n_q+1>(m-d)(q-1)-(T-d)(d-t)$ provided $q\geq d-t+1$
or $\displaystyle n_1+n_2+...+n_q+1>(m-T)(q-1)$ provided $1\leq q\leq d-t+1$.
To this end, the following proposition  below (see \cite{br} and \cite{g2}),
which follows from the properties of algebraic independent sets, will be used.

\begin{pro}
Suppose we have an algebraically independent set $A\subset\R$ and a set
$B\subset\R$ such that every element of $A$ algebraically
depends on the set $B$.
Then the cardinality of $A$ is $\leq$ the cardinality of $B$.
\end{pro}

We can assume that $\Pi^t$ and $\Pi^T$ are determined by the first $t$ and $T$
coordinates, respectively.
Let $\pi$ be the projection of $\R^m$ into the space $\R^{m-t}$ determined by
the last $m-t$ coordinates.
Then $\Pi^{d-t}=\pi(\Pi^{d})$ is a $(d-t)$-plane in $\R^{m-t}$ parallel to
the coordinate plane $\Pi^{T-t}=\pi(\Pi^T)$.
Moreover, the set of the coordinates of all points $B_{i,j}=\pi(A_{i,j})$
is algebraically independent (as a subset of the coordinates of the points
$A_{i,j}$).
Therefore, we can assume that $t=0$ by considering the space $\R^{m-t}$
determined by the last $m-t$ coordinates and the projections 
$B_{i,j}$, $\pi(\Pi^d)$, and $\pi(\Pi^T)$ into this space.

Since $Y_i\in\Pi(M_i)$, there are numbers
$\displaystyle\{\lambda_{i,j}\}_{j=1}^{n_i+1}$
such that $\displaystyle\lambda_{i,1}+...+\lambda_{i,n_i+1}=1$ and

\begin{itemize}
\item[$(1)$]
$\displaystyle Y_i=\sum_{j=1}^{n_i+1}\lambda_{i,j}A_{i,j}
=\sum_{j=1}^{n_i}\lambda_{i,j}A_{i,j}+
\big(1-\sum_{j=1}^{n_i}\lambda_{i,j}\big)A_{i,n_i+1}$ for all $i=1,...,q$.
\end{itemize}

For any $i$ at least one of the numbers
$\displaystyle\{\lambda_{i,j}\}_{j=1}^{n_i+1}$, say
$\displaystyle\lambda_{i,1}$, is different from zero
(we can even suppose that all $\displaystyle \lambda_{i,j}\neq 0$, otherwise
we reduce the sets $M_i$ by excluding the corresponding point $A_{i,j}$).
Therefore, by $(1)$, $A_{i,1}$ can be represented as a linear combination of
the points $A_{i,2},...,A_{i,n_i+1}$, $Y_i$.
Hence, $n_i$ additional numbers $\displaystyle\lambda_{i,1},...,\lambda_{i,n_i}$ are obtained.

Consequently, all coordinates of the points
$\displaystyle A_{i,j}$, $i=1,...,q$, $j=1,...,n_i+1$ are expressed in terms of
the coordinates of the points
$\displaystyle A_{i,j}$, $i=1,...,q$, $j=2,...,n_i+1$,
$\displaystyle Y_1,...,Y_q$,
and the numbers $\displaystyle\lambda_{i,j}$, $i=1,...,q$,
$\displaystyle j=1,...,n_i$.

I. Let $q\leq d+1$.
Since $\Pi^d$ is parallel to $\Pi^T$, the plane $Y_1+\Pi^T$ contains all points
$Y_2,...,Y_q$.
So,
\begin{itemize}
\item[$(2)_I$]
$\displaystyle Y_i=Y_1+\sum_{j=1}^{T}\alpha_{i,j}\bold{e}_j$, $i=2,...,q$,
\end{itemize}
where $\bold{e}_j$ denotes the $j$-th unit coordinate vector.

Therefore, any coordinate of
the points $\displaystyle A_{i,j}$, $i=1,...,q$, $j=1,...,n_i+1$, is 
algebraically dependent on the set of all coordinates of
the points $Y_1$ and $\displaystyle A_{i,j}$, $i=1,...,q$, $j=2,...,n_i+1$,
and the numbers $\displaystyle\lambda_{i,j}$, $i=1,...,q$, $j=1,...,n_i$,
$\alpha_{i,j}$, $i=2,...,q$, $j=1,...,T$.

Hence, by Proposition 2.1, $qm\leq m+n_1+...+n_q+T(q-1)$.
So,
\begin{itemize}
\item[$(3)_I$] $\displaystyle n_1+...+n_q\geq (m-T)(q-1)$.
\end{itemize}
Observe that $(3)_I$ does not change if $m$, $d$ and $T$
are replaced by $m-t$, $d-t$ and $T-t$, respectively.
Therefore, $(3)_I$ remains true for any $t $ with $0\leq t\leq d$.

\medskip
II. Let $q\geq d+1$.
Since all points $Y_1,...,Y_{q}$ are in $\Pi^d$, there are $d+1$ of them,
say $Y_1,...,Y_{d+1}$, such that each $Y_j$, $j=d+2,...,q$,
is a linear combination of the points $Y_i$, $i=1,...,d+1$.
Let us note that, such $d+1$ points exist even the linear hull of
$Y_1,...,Y_q$ is of dimension $<d$.
Therefore,
\begin{itemize}
\item[\ ]
$\displaystyle Y_j=\sum_{i=1}^{d}\beta_{j,i}Y_i+
\big(1-\sum_{i=1}^{d}\beta_{j,i}\big)Y_{d+1}$, $j=d+2,...,q$,
\end{itemize}
for some numbers $\displaystyle \{\beta_{j,i}\}$, $i=1,...,d$, $j=d+2,...,q$.

In this case equations $(2)_I$ are true for $i=2,...,d+1$.
So, any coordinate of
the points $\displaystyle A_{i,j}$, $i=1,...,q$, $j=1,...,n_i+1$,
is algebraically dependent on the set of all coordinates of 
the points $Y_1$ and $\displaystyle A_{i,j}$, $i=1,...,q$, $j=2,...,n_i+1$,
and the numbers $\displaystyle\lambda_{i,j}$, $i=1,...,q$, $j=1,...,n_i$,
$\beta_{j,i}$, $j=d+2,...,q$, $i=1,...,d$,
$\alpha_{i,j}$, $i=2,...,d+1$, $j=1,...,T$.

Hence, according to Proposition 2.1, $qm\leq m+n_1+...+n_q+(q-d-1)d+Td$,
or equivalently $n_1+...+n_q\geq (m-d)(q-1)-(T-d)d$.
Replacing in the last inequality $m$, $d$ and $T$ by $m-t$, $d-t$ and $T-t$,
respectively, we obtain that
\begin{itemize}
\item[$(3)_{II}$] $n_1+...+n_q\geq (m-d)(q-1)-(T-d)(d-t)$.
\end{itemize}
The inequalities $(3)_I$ and $(3)_{II}$ complete the proof of Theorem 1.1.

{\em Proof of Corollary $1.2$.}
We need the following considerations due to Roberts \cite{r}.
Let $\{r_i\}$ be a fixed infinite set which is algebraically independent,
i.e. every finite subset is algebraically independent.
Let $R_i=\{q+r_i:q\in \Q\}$, $\Q$ is the set of rational numbers.
Then each $R_i$ is dense in $\R$ and $R_i$'s are disjoint.
Moreover, any finite set $M$ is algebraically independent provided $M$
contains at most one point from each $R_i$, $i=1,2,...$ .

Let $\{v_i\}$ denote the vertexes of $K$ and
$R_j$ are the algebraically independent sets 
considered
above.
For each $i$ we choose a point $A_i=(A_i(1),...,A_i(m))\in\R^m$ such that
$dist\big(\theta(v_i),A_i\big)<\varepsilon $ and $A_i(s)\in R_{(i-1)m+s}$.
Then the set of all coordinates $\{A_i(k)\}$ is algebraically independent.
We define $g\colon K\to\R^m$ by $g(v_i)=A_i$ and extend $g$
linearly on every simplex of $K$.
Obviously, $g(v_{i_1})\neq g(v_{i_2})$ when $i_1\neq i_2$.
Moreover, the map $g$, restricted to any $n$-simplex of $K$ is one-to-one.
Let $\Pi^d$ be a $d$-plane in $R^m$ parallel to some coordinate planes
$\Pi^t\subset\Pi^T$ and $q$ be the number of disjoint at most $n$-dimensional
simplexes $\sigma_i=<v_{i,j}:j=1,...,n_i+1>$ from $K$ whose images 
$g(\sigma_i)=<A_{i,j}:j=1,...,n_i+1>$ under $g$ meets $\Pi^d$.
Let show that \\
$\displaystyle q\leq N_1=d+1-t+\frac{n+(n+T-m)(d-t)}{m-n-d}$
if $n\geq (m-n-T)(d-t)$.

Suppose $q\geq d+1-t$.
Since $\Pi^d$ intersects all images $g(\sigma_i)\subset
\Pi\bigl(\{A_{i,j}:j=1,...,n_i+1\}\bigr)$, $i=1,...,q$, by Theorem 1.1,
$n_1+...+n_q\geq (m-d)(q-1)-(T-d)(d-t)$.
Consequently,
\begin{itemize}
\item[(4)] $\displaystyle nq\geq (m-d)(q-1)-(T-d)(d-t)$
\end{itemize}
because $n_i\leq n$ for each $i$.
Inequality $(4)$ is equivalent to the required inequality $q\leq N_1$.

If $q\leq d+1-t$, then $q\leq N_1$ holds
because our assumption $n\geq (m-n-T)(d-t)$ implies
$d+1-t\leq N_1$.

\medskip
Now, let show that 
$\displaystyle q\leq N_2=1+\frac{n}{m-n-T}$ provided $n\leq (m-n-T)(d-t)$.

Suppose $q\leq d-t+1$.
Then, by Theorem 1.1,  $n_1+...+n_q\geq (m-T)(q-1)$.
Therefore,
\begin{itemize}
\item[(5)] $\displaystyle nq\geq (m-T)(q-1)$.
\end{itemize}
Since $(5)$ is equivalent to $q\leq N_2$, this case is completed.

Finally, let $q\geq d-t+1$. Then, as we already proved,
$q\leq N_1$. 
On the other hand, according to our assumption $n\leq (m-n-T)(d-t)$, we have 
$d-t+1\geq N_1$. Hence,
$q=d-t+1=N_1$ which implies $n=(m-n-T)(d-t)$. So,
$q=d-t+1=N_1=N_2$.


\section{Proof of Theorem 1.3}

\bigskip
First, we are going to prove Theorem 1.3 in the special case when $d=m-n$.
So, we fix spaces $X$, $Y$ and a map $f$ satisfying the hypotheses
of Theorem 1.3.
If not explicitely stated otherwise, we use the following notations
in this section:
$\varrho $ denotes the Euclidean metric in $\R^m$ and
an $\varepsilon $-disjoint set in $\R^m$ is a set which can be covered by a
family of open and disjoint subsets of $\R^m$ each of diameter $<\varepsilon $.
We say that a given set $A\subset\R^m$ is of type $(d,\varepsilon)$ if
$\Pi^d\cap A$ is $\varepsilon$-disjoint for every $d$-plane $\Pi^d$ in $\R^m$.

Let $\mathcal{H}_{\varepsilon}$, $\varepsilon >0$, denote the set of maps
$g\in C^*(X,\R^m)$ such that $ g(f^{-1}(y))$ is of type
$(m-n,\varepsilon )$ for every $y\in Y$.
Since $C^*(X,\R^m)$ equipped with the source limitation topology
has the Baire property, it suffices to show that each of the sets
$\mathcal{H}_{\varepsilon}$ is open and dense in $C^*(X,\R^m)$.
Indeed, then
$\displaystyle\mathcal{H}=\bigcap_{\displaystyle k=1}^{\infty}
\mathcal{H}_{1/k}$ would be dense and $G_{\delta}$ in $C^*(X,\R^m)$.
Moreover, if $g\in\mathcal H$ and $y\in Y$, then $ g(f^{-1}(y))\cap\Pi$
is at most 0-dimensional for every $(n-m)$-plane $\Pi\subset\R^m$.

\begin{lem}
Let $A\subset X$ be compact and $\varepsilon>0$.
Suppose the set $g_0(A)$ is of type $(d,\varepsilon)$ for some
$g_0\in C^*(X,\R^m)$ and $d$.
Then there exists a neighborhood $U$ of $A$ in $X$ and $\delta>0$ such that
$\overline{g(U)}$ is of type $(d,\varepsilon)$ provided $g\in C^*(X,\R^m)$
and $g|U$ is $\delta$-close to $g_0|U$.
\end{lem}

\begin{proof}
Assume the conclusion of the lemma is not true.
To obtain a contradiction we follow  \cite[Proof of 2.4, p. 569]{r}.
For every $i\geq 1$ take a neighborhood $U_i$ of $A$ such that
$U_i\subset g_0^{-1}(W_i)$ with $W_i$ being a
$\displaystyle 1/i$-neighborhood of $g_0(A)$.
There exist $g_i\in C^*(X,\R^m)$ and a $d$-plane $\Pi_i^d$ such that
$g_i|U_i$ is $\displaystyle 1/i$-close to $g_0|U_i$ but
$\overline{g_i(U_i)}\cap\Pi_i^d$ is not $\varepsilon$-disjoint.
Choose the points $z_i\in\overline{g_i(U_i)}\cap\Pi_i^d$ and $x_i\in U_i$
such that $\rho\big(g_i(x_i),z_i\big)\leq\displaystyle 1/i$, $i\geq 1$.
Obviously, $K=\{z_i\}_{i=1}^{\infty}\cup g_0(A)$ is a compactum
intersecting each $\Pi_i^d$.
Therefore, there exists a subsequence of $\{\Pi_i^d\}_{i=1}^{\infty}$
converging to a $d$-plane $\Pi_0^d$.
We suppose that  $\{\Pi_i^d\}_{i=1}^{\infty}$ itself converges to $\Pi_0^d$.

Let $V$ be an open subset of $\R^m$ containing $g_0(A)\cap\Pi_0^d$ and
such that it is the union of a finite disjoint and open family in $\R^m$
having  elements of diameter $<\varepsilon$.
Because each $\overline{g_i(U_i)}\cap\Pi_i^d$ is not $\varepsilon$-disjoint,
there exist points $a_i\in U_i$ and $b_i\in\overline{g_i(U_i)}\cap\Pi_i^d$
such that $V$ doesn't contain the set $\{g_i(a_i), b_i\}_{i=1}^{\infty}$.
We can also require that $\varrho\big(b_i,g_i(a_i)\big)\leq\displaystyle 1/i$
for all $i$.
This implies the existence of a point $b\in g_0(A)$ and a subsequence of
$\{b_i\}$ converging to $b$.
We still write $\displaystyle \lim b_i=b$.
Then $b\in\Pi_0^d$ because $\{\Pi_i^d\}$ converges to $\Pi_0^d$.
Hence, $b\in g_0(A)\cap\Pi_0^d\subset V$.
Consequently, $b_i\in V$ for some $i$ which contradicts the choice of $b_i$.
\end{proof}

\begin{cor}
Suppose $g_0(f^{-1}(y_0))$ is of type $(m-n,\varepsilon)$ for some $y_0\in Y$
and $g_0\in C^*(X,\R^m)$.
Then, there exists a neighborhood $V$ of $y_0$ in $Y$ and $\delta>0$
such that the set $\overline{g(f^{-1}(V))}$ is of type $(m-n,\varepsilon)$
for every $g\in C^*(X,\R^m)$ with $g|f^{-1}(V)$ being
$\delta$-close to $g_0|f^{-1}(V)$.
\end{cor}

\begin{proof}
Applying Lemma 3.1 for $f^{-1}(y_0)$, we obtain $\delta>0$ and a neighborhood
$U$ of $f^{-1}(y_0)$ such that $\overline{g(U)}$ is of type $(m-n,\varepsilon)$
provided $g\in C^*(X,\R^m)$ with $g|U$ being $\delta$-close to $g_0|U$.
Since $f$ is closed, we can find a closed neighborhood $V$ of $y_0$ in $Y$
with $f^{-1}(V)\subset U$.
Now, let $g|f^{-1}(V)$ be $\delta$-close to $g_0|f^{-1}(V)$
for some $g\in C^*(X,\R^m)$.
Extend $g|f^{-1}(V)$ to a map $h\in C^*(X,\R^m)$ such that
$h|U$ remains $\delta$-close to $g_0|U$.
Then, according to the choice of $U$ and $\delta$,
$\overline{h(U)}$ is of type $(m-n,\varepsilon)$.
Finally, since $\overline{h(f^{-1}(V))}\subset \overline{h(U)}$ and
$h|f^{-1}(V)=g|f^{-1}(V)$, we are done.
\end{proof}

\begin{pro}
Any $\mathcal{H}_{\varepsilon}$ is open in $C^*(X,\R^m)$.
\end{pro}

\begin{proof}
We fix $g_0\in \mathcal{H}_{\varepsilon}$.
By Corollary 3.2, for every $y\in Y$ there exists a neighborhood $V_y$ of $y$
in $Y$ and $\delta_y >0$ such that if $g\in C^*(X,\R^m)$ and
$\varrho\big(g(x),g_0(x)\big)\leq\delta_y$ for all $x\in f^{-1}(V_y)$,
then $\overline{g(f^{-1}(V_y))}$ is of type $(m-n,\varepsilon)$.
Take a locally finite open cover $\omega$ of $Y$ refining $\{V_y: y\in Y\}$
and for each $W\in\omega $ fix $y(W)\in Y$ such that $W\subset V_{y(W)}$.
Define the set-valued map $\phi\colon Y\to\ (0,\infty)$ by
$\phi(y)=\cup\{(0,\delta_{y(W)}]: y\in W\}$.
Obviously, $\phi $ is convex-valued and lower semi-continuous.
By \cite[Theorem 6.2]{rs}, $\phi$ admits a continuous selection
$\beta\colon Y\to\ (0,\infty)$ and let $\alpha =\beta \circ f$.
It suffices to show that  $g\in C^*(X,\R^m)$ and
$\varrho\big(g_0(x),g(x)\big)<\alpha (x)$ for every $x\in X$ imply
$g\in \mathcal{H}_{\varepsilon}$.
To this end, let $y\in Y$ and select $W\in\omega$ containing $y$ such that
$\alpha (x)\leq\delta_{y(W)}$ for every $x\in f^{-1}(y)$.
Take a function $h_y\in C^*(X,\R^m)$ coinciding with $g$
on the set $f^{-1}(y)$ and satisfying the inequality
$\varrho\big(h_y(x),g_0(x)\big)\leq\delta_{y(W)}$ for all $x\in X$.
Then, according to the choice of $V_{y(W)}$, the set
$\overline{h_y(f^{-1}(V_{y(W)}))}$ is of type $(m-n,\varepsilon)$,
so is $\overline{g(f^{-1}(y))}$.
Hence, any $g\in C^*(X,\R^m)$ which is $\alpha$-close to $g_0$ belongs to
$\mathcal{H}_{\varepsilon}$.
Therefore, $\mathcal{H}_{\varepsilon}$ is open in $C^*(X,\R^m)$.
\end{proof}

For every space $M$ and $\varepsilon>0$ let
$\displaystyle C_{(n,\varepsilon)}(M,\R^m)$ denote the set of all
$g\in C^*(M,\R^m)$ such that $g(M)$ is of type
$(m-n,\varepsilon)$ provided $m\geq n+1$.

\begin{lem}
Let $M$ be an $n$-dimensional compactum and $m\geq n+1$.
Then $\displaystyle C_{(n,\varepsilon)}(M,\R^m)$ is dense in $C(M,\R^m)$
for every $\varepsilon>0$.
\end{lem}

\begin{proof}
Let $g_0\in C(M,\R^m)$ and $\delta>0$.
Representing $g_0$ as the composition of two maps $q_1\colon M\to Z$ and
$q_2\colon Z\to\R^m$, where $Z$ is a metrizable compactum of dimension
$\leq n$, and considering $Z$ and $q_2$ instead of $M$ and $g_0$,
we reduce our proof to the case $M$ is a metrizable compactum.
Take a positive number $\eta$ satisfying the following conditions:

\begin{itemize}
\item[(6)] \hbox{}~~~~~~$\displaystyle 5\eta<\delta/2$ and
$\displaystyle 9\eta (r+1)<\varepsilon$, where $\displaystyle r=n(m+1-n)$.
\end{itemize}

\medskip\noindent
Since $\dim M\leq n$, by a standard procedure we can find a finite
$n$-dimensional complex $K$ and maps $h\colon M\to K$ and
$\theta \colon K\to \R^m $ such that $\theta \circ h$ is
$\displaystyle\delta/2$-close to $g_0$.
Moreover, we can suppose that:

\begin{itemize}
\item[(7)] $diam\bigl(\theta (\sigma)\bigr)<\eta$ for every
simplex $\sigma\in K$.
\end{itemize}

\noindent
It suffices to find a map $g\colon K\to\R^m$ which is
$\displaystyle\delta/2$-close to the map $\theta $ and
$g\circ h\in \displaystyle C_{(n,\varepsilon)}(M,\R^m)$.
To this end, we apply Corollary 1.2
(with $d=m-n$, $t=0$, $T=m$, $\varepsilon =\eta $ and $n$ replaced by $n-1$)
to obtain a semi-linear map $g\colon K\to\R^m$ such that
$\varrho\bigl(g(v),\theta (v)\bigr)\leq \eta $ for all vertexes $v$ of $K$ and,
for each $(m-n)$-plane $\Pi\subset\R^m$, the number $q$ of disjoint at most
$(n-1)$-dimensional simplexes of $K$ whose images under $g$ intersect $\Pi$
is $\leq r=n(m+1-n)$.
We can choose $g$ such that, in addition to the above requirements, we
also have $g(v_i)\neq g(v_j)$ for any different vertices $v_i$ and $v_j$ of $K$.

\noindent
Let $v_i$ and $v_j$ be two vertices of the same simplex $\triangle\in K$.
Then, according to (7) and the choice of $g$, we have

\medskip\noindent
$\varrho\bigl(g(v_i),g(v_j)\bigr)\leq \varrho\bigl(g(v_i),\theta (v_i)\bigr)+
\varrho\bigl(\theta (v_i), \theta (v_j)\bigr)+
\varrho\bigl(\theta (v_j),g(v_j)\bigr)<3\eta $.

\medskip\noindent
Consequently

\begin{itemize}
\item[(8)] $g(\triangle)$ is of diameter $<3\eta $ for any simplex
$\triangle\in K$.
\end{itemize}

\noindent
Item (8) implies that $\varrho\bigl(g(y),\theta(y)\bigr)<5\eta$
for all $y\in K$.
Hence, by (6), $g$ and $\theta$ are $\displaystyle\delta/2$-close to each other.

\medskip\noindent
It remains only to show that
$g\circ h\in \displaystyle C_{(n,\varepsilon)}(M,\R^m)$, or equivalently,
$g(K)$ is of type $(m-n,\varepsilon)$.
To this end we use an idea of \cite[proof of 2.3, p. 568]{r}.
We fix a $(m-n)$-plane $\Pi\subset\R^m$.
It is enough to prove that each component of $g(K)\cap\Pi$ is of diameter
$\leq 9(r+1)\eta$ because, by (6), $9(r+1)\eta<\varepsilon$.
Suppose $\varrho(a,b)>9(r+1)\eta$ for some component $P$ of $g(K)\cap\Pi$
and some points $a, b\in P$.
There is an arc $ab$ in $P$.
Because of (8), every subarc of $ab$ of diameter $\geq 3\eta$ must contain
at least one point on the boundary of a simplex of $g(K)$,
hence on a simplex of $g(K)$ of dimension $\leq n-1$.
Next, we take points $a_i\in ab$, $i=1,...,r+1$ such that: \\
$3(3i-2)\eta< \varrho(a,a_i)\leq 3(3i-1)\eta$, $i=1,...,r+1$, and each
$a_i$ belongs to a simplex $g(\sigma_i)\in g(K)$ with $\sigma_i$ being of dimension $\leq n-1$. \\
Then, $\varrho(a_i,a_j)>6\eta$ for $i\neq j$ which, according to (8),
implies $g(\sigma_i)\cap g(\sigma_j)=\emptyset$.
Therefore, we obtained $r+1$ disjoint simplexes $\sigma_i\in K$, $i=1,...,r+1$,
each of them having dimension $\leq n-1$ whose images under $g$ meets $\Pi$.
This is a contradiction because, according to the choice of $g$,
the number of pairwise disjoint at most $(n-1)$-dimensional simplexes of $K$
whose images under $g$ intersect $\Pi$ is $\leq n(m-n+1)=r$.
\end{proof}

Below we consider the set-valued map
$\displaystyle\psi_{\varepsilon}\colon Y\to C^*(X,\R^m)$  defined by
$\psi_{\varepsilon}(y)=C^*(X,\R^m)\backslash\mathcal{H}_{\varepsilon}(y)$,
where $\mathcal {H}_{\varepsilon}(y)$ denotes the set of all
$g\in C^*(X,\R^m)$ such that $g(f^{-1}(y))$ is of type $(m-n,\varepsilon)$.

Next proposition, combined with Proposition 3.3, will complete
the proof of Theorem 1.3 for $d=m-n$.

\begin{pro}
Any $\mathcal{H}_{\varepsilon}$ is dense in $C^*(X,\R^m)$.
\end{pro}

\begin{proof}
We first show that the graph $G$ of $\psi_{\varepsilon}$ is closed in
$Y\times C^*(X,\R^m)$ provided $C^*(X,\R^m)$ is equipped with
the uniform convergence topology generated by the metric $\varrho $.
Let $(y_0,g_0)\in (Y\times C^*(X,\R^m))\backslash G$.
Then $g_0\not\in\psi_{\varepsilon}(y_0)$, so
$g_0\in\mathcal{H}_{\varepsilon}(y_0)$.
By Corollary 3.2, there exists $\delta>0$ and a neighborhood $V$ of $y_0$
in $Y$ such that $\overline{g(f^{-1}(V))}$ is of type $(m-n,\varepsilon)$
for every $g\in C^*(X,\R^m)$ with $g|f^{-1}(V)$ being $\delta$-close to
$g_0|f^{-1}(V)$.
Let $W\subset C^*(X,\R^m)$ be the set of all $g$ which are $\delta$-close
to $g_0$.
Obviously, $W$ is a neighborhood of $g_0$ in $C^*(X,\R^m)$ with respect to
the uniform convergence topology, so $V\times W$ is a neighborhood of
$(y_0,g_0)$ in $Y\times C^*(X,\R^m)$ which doesn't meet $G$.
Therefore, $G$ is closed in $Y\times C^*(X,\R^m)$.

{\em Claim.
$\overline{B}(g_0,\alpha)\backslash\psi_{\varepsilon}(y)\neq\emptyset$
for every $y\in Y$, $\alpha\colon X\to (0,\infty)$ and $g_0\in C^*(X,\R^m)$.}

We need to show that
$\overline{B}(g_0,\alpha)\cap\mathcal{H}_{\varepsilon}(y)\neq\emptyset$
for fixed $y\in Y$, $\alpha\in C(X,(0,\infty))$ and $g_0\in C^*(X,\R^m)$.
Let $\delta >0$ be the minimal value of $\alpha $ on $f^{-1}(y)$.
Since $\dim f^{-1}(y)\leq n$, by Lemma 3.4, there exists
$h\in C(f^{-1}(y),\R^m)$ which is $\delta$-close to $g_0|f^{-1}(y)$ and
such that $h(f^{-1}(y))$ is of type $(m-n,\varepsilon)$.
Obviously, $\varrho\big(h(x),g_0(x)\big)\leq\alpha(x)$ for all $x\in f^{-1}(y)$.
So, every extension $g\in C^*(X,\R^m)$ of $h$ would be in
$\mathcal{H}_{\varepsilon}(y)$.
Hence, the proof is reduced to find such an extension of $h$ which is also in
$\overline{B}(g_0,\alpha)$.
To this end, we define the set-valued map $\Phi\colon X\to\R^m$ by
$\Phi(x)=h(x)$ if $x\in f^{-1}(y)$ and
$\Phi(x)=B\big(g_0(x),\alpha(x)\big)$ if $x\not\in f^{-1}(y)$.
Here $B\big(g_0(x),\alpha(x)\big)$ denotes the closed ball in $\R^m$
having center $g_0(x)$ and radius $\alpha(x)$.
This map is lower semi-continuous with closed and convex values in $\R^m$.
Hence, by Michael's convex-valued selection theorem, $\Phi$ admits a
continuous selection $g$ which extends $h$.
Moreover, we can assume that $\alpha$ is a bounded function, which implies
that $g$ is also bounded.
Therefore, $g\in\overline{B}(g_0,\alpha)$, which completes the proof
of the claim.

\smallskip
We are now in a position to finish the proof of our proposition.
Fix $g_0\in C^*(X,\R^m)$ and $\alpha\in C^*(X,(0,1))$, and consider the
constant set-valued map $\varphi\colon Y\to C^*(X,\R^m)$,
$\varphi(y)=\overline{B}(g_0,\alpha)$, where $C^*(X,\R^m)$ is equipped with
the uniform convergence topology.
Now, we need the following theorem of Michael \cite[Theorem 5.3]{m}: \\
Let $Y$ be paracompact with $\dim Y=0$, $M$ a completely metrizable,
and $\phi\colon Y\to M$ a lower semi-continuous closed-valued map.
If $\psi\colon Y\to M$ is a set-valued map with a closed graph
such that $\phi (y)\backslash\psi(y)\neq\emptyset$ for every $y\in Y$,
then $\phi$ has a selection avoiding $\psi$. \\
In our case $\varphi$ and $\psi_{\varepsilon}$ satisfy the hypotheses of the
Michael's theorem, so there exists a map $\theta\colon Y\to C^*(X,\R^m)$ such
that $\theta (y)\in \overline{B}(g_0,\alpha)\backslash\psi_{\varepsilon}(y)$
for all $y\in Y$.
We define the map $g\in C^*(X,\R^m)$, $g(x)=\theta(f(x))$, $x\in X$.
Then $g\in\overline{B}(g_0,\alpha)\cap\mathcal{H}_{\varepsilon}(y)$
for every $y\in Y$.
Consequently, $g\in\overline{B}(g_0,\alpha)\cap\mathcal{H}_{\varepsilon}$
and we are done.
\end{proof}

As we already mentioned, it follows from Proposition 3.3 and Proposition 3.5
that the set $\mathcal{H}$ is dense and $G_{\delta}$ in $C^*(X,\R^m)$.
This provides the proof of Theorem 1.3 in the special case $d=m-n$.
Let us prove the general case of Theorem 1.3.
We are going to show that every $g\in\mathcal H$ satisfies the theorem, i.e. it 
satisfies the following condition:
$A=g(f^{-1}(y))\cap\Pi^d$ is at most $(n+d-m)$-dimensional for every
$y\in Y$ and every $d$-plane $\Pi^d\subset\R^m$ with $m-n+1\leq d\leq m$.
Fix some $(m-n)$-plane $\Pi^{m-n}\subset \Pi^d$ and consider the orthogonal
projection $p$ of $\Pi^d$ onto the $(n+d-m)$-plane
$\Pi^{n+d-m}\subset \Pi^d$ which is the orthogonal complement of 
$\Pi^{m-n}$ in $\Pi^d$.
Then the compactum $B=p\bigl(g(f^{-1}(y))\cap\Pi^d\bigr)\subset \Pi^{n+d-m}$
is of dimension $\leq n+d-m$.
Obviously, any fiber of the mapping $p|_A\colon A\to B$ is the intersection of
$g(f^{-1}(y))$ and some $(m-n)$-plane. So, according to
the choice of $g$, $p|_A$ has zero-dimensional fibers.
Therefore, by Hurewicz's theorem about perfect zero-dimensional maps, 
$\dim A\leq \dim B\leq n+d-m$.

\section{Proof of Theorem 1.5 and Theorem 1.7}

\bigskip\noindent
{\em Proof of Theorem $1.5$.}
It suffices to show that, for any integers $d,t,T$ with
$0\leq t\leq d\leq T\leq m$ and $d\leq m-n-1$, and any coordinate planes
$\Pi^t\subset\Pi^T$ in $\R^m$, $C^*(X,\R^m)$ contains a dense
$G_{\delta}$-subset of maps $g$ such that $f^{-1}(y)\cap g^{-1}(\Pi^d)$
has at most $q$ points for every $y\in Y$ and every $d$-plane
$\Pi^d\subset\R^m$ parallel to $\Pi^t\subset\Pi^T$,
where $q$ is the integer part of
$\displaystyle N_1=d+1-t+\frac{n+(n+T-m)(d-t)}{m-n-d}$ if $n\geq (m-n-T)(d-t)$
and the integer part of $\displaystyle N_2=1+\frac{n}{m-n-T}$ otherwise.
So, we fix integers $n,d,t,T$ satisfying the above inequalities and coordinate
planes $\Pi^t\subset\Pi^T\subset\R^m$.
Observe that always $q\geq 1$, and let $\mathcal{P}$ be the set of all
$d$-planes in $\R^m$ parallel to $\Pi^t\subset\Pi^T$.
We now define the following augmented notion which is used in this section:
A subset $A$ of an arbitrary metric space $M$ is said to be of {\em cotype}
$(q,\varepsilon, g)$, where $\varepsilon >0$ and $g\in C^*(M,\R^m)$,
if for any $\Pi^d\in\mathcal{P}$, the set $A\cap g^{-1}(\Pi^d)$ can be covered
by at most $q$ disjoint open sets in $M$ each of diameter $\leq\varepsilon$.
Let $\mathcal{K}_{\varepsilon}$, $\varepsilon>0$, denote the set of all
maps $g\in C^*(X,\R^m)$ such that $f^{-1}(y)$ is of cotype $(q,\varepsilon,g)$
for every $y\in Y$.
The proof is reduced to show that each $\mathcal{K}_{\varepsilon}$
is open and dense in $C^*(X,\R^m)$.

\begin{lem}
Let $A$ be a compact subset of $X$ which is of cotype $(q,\varepsilon, g_0)$
for some $g_0\in C^*(X,\R^m)$ and $\varepsilon>0$.
Then there exists a neighborhood $U$ of $A$ in $X$ and $\delta>0$
with $U$ being of cotype $(q,\varepsilon,g)$ for every $g\in C^*(X,\R^m)$
such that $g|U$ is $\delta$-close to $g_0|U$.
\end{lem}

\begin{proof}
Suppose the lemma is not true.
For every $i\geq 1$ take a $\displaystyle 1/i$-neighborhood $U_i$ of $A$
such that $g_0(U_i)$ is within $\displaystyle 1/i$-neighborhood of $g_0(A)$
in $\R^m$.
Then, there exist $g_i\in C^*(X,\R^m)$ and a $d$-plane $\Pi_i^d\in\mathcal{P}$
such that $g_i|U_i$ is $\displaystyle 1/i$-close to $g_0|U_i$, but
$g_i^{-1}(\Pi_i^d)\cap U_i$ is not covered by any family of $\leq q$
open disjoint sets in $X$ with diameters $\leq\varepsilon$.
As in the proof of Lemma 3.1, passing to subsequences, we can suppose that
$\{\Pi_i^d\}_{i=1}^{\infty}$ converges to a $d$-plane $\Pi_0^d$.
Since all $\Pi_i^d$ are from $\mathcal{P}$
(i.e., parallel to $\Pi^t\subset\Pi^T$), so is $\Pi_0^d$.
Hence, $A\cap g_0^{-1}(\Pi_0^d)$ is covered by a disjoint family
$\{V_j\}$ of $\leq q$ open sets in $X$ with diameters $\leq\varepsilon$.
Take points $x_i\in g_i^{-1}(\Pi_i^d)\cap U_i\backslash V$ and $y_i\in A$
such that $\displaystyle dist(x_i,y_i)\le 1/i$, $i=1,2,...$,
where $V=\bigcup V_j$.
We can assume that the sequence ${x_i}$ converges to some $x_0\in A$
(recall that $A$ is compact).
Then $\{g_i(x_i)\}$ converges to $g_0(x_0)$ and $g_0(x_0)\in g_0(A)\cap\Pi_0^d$.
Therefore, $x_0\in A\cap g_0^{-1}(\Pi_0^d)\subset V$.
So, $x_i\in V$ for almost all $i$ which is a contradiction.
\end{proof}

The proof of next corollary is similar to that one of Corollary 3.2.
\begin{cor}
Suppose $f^{-1}(y_0)$ is of cotype $(q,\varepsilon, g_0)$ for some $y_0\in Y$
and $g_0\in C^*(X,\R^m)$.
Then, there exists a neighborhood $V$ of $y_0$ in $Y$ and $\delta>0$
such that the set $f^{-1}(V)$ is of cotype $(q,\varepsilon, g)$ for every
$g\in C^*(X,\R^m)$ with $g|f^{-1}(V)$ being $\delta$-close to  $g_0|f^{-1}(V)$.
\end{cor}

\begin{pro}
Any $\mathcal{K}_{\varepsilon}$ is open in $C^*(X,\R^m)$.
\end{pro}

\begin{proof}
The proof follows the arguments from the proof of Proposition 3.3,
but now we apply Corollary 4.2 instead of Corollary 3.2.
\end{proof}

\begin{lem}
Let $M$ be a metrizable at most $n$-dimensional compactum.
Then the set $\mathcal{K}_0(M,\R^m)$ of all $g\in C(M,\R^m)$
such that $g^{-1}(\Pi^d)$ consists of at most $q$ points for every
$\Pi^d\in\mathcal{P}$, is dense in $C(M,\R^m)$.
\end{lem}

\begin{proof}
Let $\Omega$ be the collection of all disjoint families
$\{\overline{V}_1,\overline{V}_2,...,\overline{V}_{q+1}\}$ of $q+1$ elements
such that each $V_j$ belongs to a fixed base for $M$.
Let also

\smallskip\noindent
$C_{\Gamma}=\{g\in C(M,\R^m): g^{-1}(\Pi^d)
\hbox{}~\mbox{meets at most}\hbox{}~q\hbox{}~\mbox{elements of}\hbox{}~
\Gamma\hbox{}~\mbox{for every}\hbox{}~\Pi^d\in\mathcal{P}\}$, $\Gamma\in\Omega$.

Obviously, $\Omega$ is countable and $\mathcal{K}_0(M,\R^m)$ is
the intersection of all sets $C_{\Gamma}$, $\Gamma\in\Omega$.
Therefore, our proof is reduced to show that each
$C_{\Gamma}$ is dense and open in $C(M,\R^m)$.

\medskip\noindent
\hspace{0.3cm}{\em Claim $1$.} Every $C_{\Gamma}$ is open in $C(M,\R^m)$.

\medskip\noindent
We fix $\Gamma\in\Omega$ and $g_0\in C_{\Gamma}$.
Suppose, for every $i$ there exists $g_i\not\in C_{\Gamma}$
with $g_i$ being $1/i$-close to $g_0$.
So, we can find $\Pi_i^d\in\mathcal{P}$ such that $g_i^{-1}(\Pi_i^d)$
meets every element of $\Gamma$.
As in Lemma 3.1, we can suppose that the sequence $\{\Pi_i^d\}$
converges to some $\Pi_0^d\in\mathcal{P}$.
Then $g_0^{-1}(\Pi_0^d)$ intersects at most $q$ elements of $\Gamma$,
let say the first $q$.
Now, for every $i$, choose a point
$x_i\in g_i^{-1}(\Pi_i^d)\cap\overline{V}_{q+1}$ and, since $M$ is compact,
we can assume that the sequence $\{x_i\}$ converges to some
$x_0\in\overline{V}_{q+1}$.
Then, $\{g_i(x_i)\}$ converges to $g_0(x_0)\in\Pi_0^d$.
So, $x_0\in g_0^{-1}(\Pi_0^d)\cap\overline{V}_{q+1}$ which is a contradiction.

\medskip\noindent
\hspace{0.3cm}{\em Claim $2$.} Every $C_{\Gamma}$ is dense in $C(M,\R^m)$.
\medskip\noindent

Let $\Gamma=\{\overline{V}_1,\overline{V}_2,...,\overline{V}_{q+1}\}$,
$g_0\in C_{\Gamma}$ and $\delta>0$.
There exist an open cover $\omega$ of $M$ with $\mesh(\omega)\leq r/3$,
where $r=\min\{dist(\overline{V}_i,\overline{V}_j): i\neq j\}$,
and a semi-linear map $h\colon L \to\R^m$ such that
$g=h\circ\pi$ is $\delta$-close to $g_0$.
Here $L$ is the polyhedron underlying the nerve of $\omega$ and
$\pi\colon M\to L$ is the canonical map.
According to Corollary 1.2, we can assume that for any $\Pi^d\in\mathcal{P}$,
the number of pairwise disjoint, at most $n$-dimensional simplexes
$\sigma\in L$ with $h(\sigma)$ meeting $\Pi^d$, is $\leq q$.
We can also assume that $\omega$ is of order $\leq n+1$,
so $L$ is at most $n$-dimensional.
If there exists a plane $\Pi^*\in\mathcal{P}$ with $g^{-1}(\Pi^*)$
intersecting every $\overline{V}_i$, choose
$x_i\in g^{-1}(\Pi^*)\cap\overline{V}_i$, and let $\omega_i$ be the family of
those elements of $\omega$ containing the point $x_i$.
Then each family $\omega_i$, $i=1,...,q+1$, generates a simplex $\sigma_i\in L$
of dimension $\leq n$ such that $h(x_i)\in h(\sigma_i)\cap \Pi^*$
and $\sigma_i\cap\sigma_j=\emptyset$ for $i\neq j$.
This contradicts the choice of $h$.
\end{proof}

Since $\mathcal{K}_0(M,\R^m)\subset\mathcal{K}_{\varepsilon}(M,\R^m)$ for
every $\varepsilon>0$, where $\mathcal{K}_{\varepsilon}(M,\R^m)$ is the set
of all maps $g\in C^*(M,\R^m)$ with $M$ being of cotype $(q,\varepsilon,g)$,
Lemma 4.4 implies the following corollary:

\begin{cor}
Let $M$ be a metrizable compactum of dimension $\leq n$.
Then any $\mathcal{K}_{\varepsilon}(M,\R^m)$ is dense in $C(M,\R^m)$.
\end{cor}

Finally, next proposition, combined with Proposition 4.3,
completes the proof of Theorem 1.5.

\begin{pro}
Any $\mathcal{K}_{\varepsilon}$ is dense in $C^*(X,\R^m)$.
\end{pro}

\begin{proof}
Let $\displaystyle\psi_{\varepsilon}\colon Y\to C^*(X,\R^m)$ be the set-valued
map $\psi_{\varepsilon}(y)=C^*(X,\R^m)\backslash\mathcal{K}_{\varepsilon}(y)$,
where $\mathcal {K}_{\varepsilon}(y)$ denotes the set of all
$g\in C^*(X,\R^m)$ such that $f^{-1}(y)$ is of cotype $(q,\varepsilon,g)$.
Further, the proof follows the arguments from the proof of Proposition 3.5.
To show that the graph of $\psi_{\varepsilon}$ is closed,
we use now Corollary 4.2 instead of Corollary 3.2.
Also, the application of Lemma 3.4 should be replaced by that one of
Corollary 4.5 in the proof of the claim.
\end{proof}

\bigskip\noindent
{\em Proof of Theorem $1.7$.}
We are going to prove that $C^*(X,l_2)$ contains a dense $G_{\delta}$-subset
of maps $g$ such that for every $y\in Y$ and every $\Pi^d\in \mathcal{P}(d,r)$,
$f^{-1}(y)\cap g^{-1}(\Pi^d)$ has at most
$\displaystyle d+1-r$ points if $r\leq d$ and at most one point otherwise.
For given integers $d,r$ let $\mathcal{Q}$ be the set of all $d$-planes
parallel to the $r$-plane of first $r$ coordinates in $l_2$.
As in Theorem 1.5, we introduce the notion of a $(q,\varepsilon, g)$-cotype
subset of a metric space $M$ considering now planes $\Pi^d\in\mathcal{Q}$
and taking $q=d+1-r$ in case $d\geq r$ and $q=1$ otherwise.
Let $\mathcal{F}_{\varepsilon }$, $\varepsilon>0$, denote the set of all
maps $g\in C^*(X,l_2)$ such that $f^{-1}(y)$ is of cotype $(q,\varepsilon,g)$
for every $y\in Y$.
It suffices to
show that each $\mathcal{F}_{\varepsilon }$ is open and dense in $C^*(X,l_2)$.
Following the proof of Theorem 1.5, we can show that all
$\mathcal{F}_{\varepsilon }$ are open in $C^*(X,l_2)$.
To prove density of $\mathcal{F}_{\varepsilon }$, we use the following lemma.

\begin{lem}
Let $M$ be a metrizable compactum and $\varepsilon>0$.
Then the set $\mathcal{F}_{\varepsilon}(M,l_2)$ of all $g\in C(M,l_2)$
such that $M$ is of cotype $(q,\varepsilon,g)$, is dense in $C(M,l_2)$.
\end{lem}

\begin{proof}
We take $g_0\in C(M,l_2)$ and $\lambda>0$.
Then, there exist maps $\varphi\colon M\to K$ and $h\colon K\to l_2$ such that
$K$ is a finite complex and $h\circ\varphi$ is $\lambda/2$-close to $g_0$.
Moreover, we can assume that each fiber of $\varphi$ has diameter
$<\varepsilon $.
For every $A\subset\N$ we identify $\R^A$ with the subspace
$\{y\in l_2: y_i=0\hbox{}~\mbox{for all}\hbox{}~i\not\in A\}$ of $l_2$,
and denote by $\pi_A$ the canonical projection $\pi_A\colon l_2\to\R^{A}$.
Let $\dim K=n$ and $A\subset\N$ be any finite set satisfying
the following conditions: \\
(i) $\{1,...,r\}\subset A$;  \\
(ii) $|A|\geq r+d+2n+n\cdot |d-r|+1$, where $|A|$ is the cardinality of $A$. \\
Since every projection $\pi_A$ is open, we can show that the map
$\Lambda_A\colon C(K,l_2)\to C(K,\R^A)$, $\Lambda_A(g)=\pi_A\circ g$,
is also open.

\smallskip
Suppose first that $r\leq d$.
Observe that, if $T=m$, then $(m-n-T)(d-r)=-n(d-r)\leq 0$,
so $n\geq (m-n-T)(d-r)$.
Hence, applying Theorem 1.5 with $m=T=|A|$ and $t=r$,
there is a dense $G_{\delta}$ subset $\mathcal{F}_A$ of $C(K,\R^A)$
consisting of maps $g$ such that $g^{-1}(\Pi^d)$ contains at most
$\displaystyle 1+d-r+\frac{n+n(d-r)}{|A|-n-d}$ points for every $d$-plane
$\Pi^d\subset\R^A$ parallel to the $r$-plane of first $r$ coordinates in $\R^A$.
Since $|A|\geq r+d+2n+n\cdot (d-r)+1$, any preimage $g^{-1}(\Pi^d)$, $g\in\mathcal{F}_A$, 
contains at most $\displaystyle 1+d-r$ points.
The set $\Lambda_A^{-1}\big(\mathcal{F}_A\big)$ is also dense and
$G_{\delta}$ in $C(K,l_2)$.
Therefore, there exists $g\in \Lambda_A^{-1}\big(\mathcal{F}_A\big)$
which is $\lambda/2$-close to $h$.
Then $\overline{g}=g\circ\varphi$ is $\lambda$-close to $g_0$.
It remains only to show that $\overline{g}\in\mathcal{F}_{\varepsilon }(M,l_2)$.
To this end, take $\Pi^d\in\mathcal{Q}$.
Since $\Pi_A^d=\pi_A(\Pi^d)$ is a $d$-plane in $\R^A$ parallel to the
$r$-plane of first $r$ coordinates in $\R^A$, 
$g^{-1}\big(\pi_A^{-1}(\Pi_A^d)\big)$ contains $\displaystyle\leq 1+d-r$ points.
Formally the plane $\Pi_A^d$ can has dimension $<d$,
but in this case the estimation on the cardinality of preimage is even better.
The inclusion $\Pi^d\subset \pi_A^{-1}(\Pi_A^d)$ implies that $g^{-1}(\Pi^d)$
contains $\leq 1+d-r$ points.
Hence, $(\overline{g})^{-1}(\Pi^d)=\varphi^{-1}\big(g^{-1}(\Pi^d)\big)$
consists of $\leq 1+d-r$ fibers of $\varphi $.
Since any fiber of $\varphi $ is of diameter $<\varepsilon $,
$\overline{g}\in\mathcal{F}_{\varepsilon }(M,l_2)$.

\smallskip
Suppose $d\leq r$.
We again apply Theorem 1.5 to $A$, but now for $m=|A|$, $t=0$ and $T=r$.
Obviously, in this case $(m-n-T)(d-t)>n$.
Therefore, there is a dense $G_{\delta}$ subset $\mathcal{F}_A$ of $C(K,\R^A)$
consisting of maps $g$ such that $g^{-1}(\Pi^d)$ contains at most one point
for every $d$-plane in $\R^A$ parallel to the $r$-plane of first $r$
coordinates in $\R^A$.
As before, we take $g\in \Lambda_A^{-1}\big(\mathcal{F}_A\big)$
which is $\lambda/2$-close to $h$ and show that
$\overline{g}=g\circ\varphi\in\mathcal{F}_{\varepsilon }(M,l_2)$.
\end{proof}

\begin{pro}
Any $\mathcal{F}_{\varepsilon}$ is dense in $C^*(X,l_2)$.
\end{pro}

\begin{proof}
We follow the proof of Proposition 3.5, describing the necessary changes.
In our situation $\displaystyle\psi_{\varepsilon}\colon Y\to C^*(X,l_2)$
is defined by
$\psi_{\varepsilon}(y)=C^*(X,l_2)\backslash\mathcal{F}_{\varepsilon }(y)$,
where $\mathcal {F}_{\varepsilon }(y)$ denotes the set of all
$g\in C^*(X,l_2)$ such that $f^{-1}(y)$ is of cotype $(q,\varepsilon,g)$
with $q=1+d-r$ if $r\leq d$ and $q=1$ otherwise.
To show that $\displaystyle\psi_{\varepsilon}$ has a closed graph,
we apply the analogue of Corollary 4.2 for the space $C^*(X,l_2)$
instead of Corollary 3.2.
We also need the following claim.

\medskip\noindent
\hspace{0.3cm}{\em Claim.}
{\em Let $g_0\in C^*(X,l_2)$, $\alpha\colon X\to (0,\infty)$ and $y\in Y$.
Then, for any $\varepsilon>0$,
$\psi_{\varepsilon}(y)\cap\overline{B}(g_0,\alpha)$ is a $Z$-set in
$\overline{B}(g_0,\alpha)$ provided $\overline{B}(g_0,\alpha)$ is considered
as a subset of $C^*(X,l_2)$ with the uniform convergence topology.}

\medskip\noindent
The proof of this claim partly follows the arguments
from the proof of \cite[Lemma 8]{tv}.
We need to show that every map $h\colon Q\to\overline{B}(g_0,\alpha)$,
where $Q$ is the Hilbert cube, can be approximated by a map
$h_1\colon Q\to \overline{B}(g_0,\alpha)$ avoiding the set
$\psi_{\varepsilon}(y)\cap\overline{B}(g_0,\alpha)$.
So, we fix such a map $h$ and $\eta>0$.
Then $h$ generates the map $u\colon Q\times X\to l_2$, $u(z,x)=h(z)(x)$,
such that $dist(u(z,x),g_0(x))\leq\alpha(x)$ for any $(z,x)\in Q\times X$.
Choose $\lambda\in (0,1)$ with $\lambda\sup\{\alpha(x):x\in f^{-1}(y)\}<\eta/2$
and define $u_1\in C(Q\times f^{-1}(y),l_2)$ by
$u_1(z,x)=(1-\lambda)u(z,x)+\lambda g_0(x)$.
For every $(z,x)\in Q\times f^{-1}(y)$ we have

\hspace{0.4cm} $dist(u_1(z,x),g_0(x))<\alpha(x)$ and
$dist(u_1(z,x),u(z,x))<\eta/2$.

\smallskip\noindent
Let $\delta<\inf\{\eta/2,\alpha(x)-dist(u_1(z,x),g_0(x)):
(z,x)\in Q\times f^{-1}(y)\}$.
By Lemma 4.7, there exists a map $u_2\colon Q\times f^{-1}(y)\to l_2$
such that $Q\times f^{-1}(y)$ is of cotype $(q,\varepsilon,u_2)$ and
$dist(u_1(z,x),u_2(z,x))<\delta$ for all $(z,x)\in Q\times f^{-1}(y)$.
Then

\hspace{0.4cm} $dist(u_2(z,x),g_0(x))<\alpha(x)$ and
$dist(u_2(z,x),u(z,x))<\eta$

\smallskip\noindent
provided $(z,x)\in Q\times f^{-1}(y)$.
The equality $h_2(z)(x)=u_2(z,x)$ determines a map
$h_2\colon Q\to C(f^{-1}(y),l_2)$ such that
$f^{-1}(y)$ is of cotype $(q,\varepsilon,h_2(z))$ for every $z\in Q$.
One can show that the projection map
$pr\colon\overline{B}(g_0,\alpha)\to C(f^{-1}(y),l_2)$, $pr(g)=g|f^{-1}(y)$,
is open with respect to the uniform convergence topology and
$pr\big(\overline{B}(g_0,\alpha)\big)$ contains $h_2(Q)$.
This yields that $h_2$ can be lifted to a map
$w\colon Q\to\overline{B}(g_0,\alpha)$ such that $w$ is $\eta$-close to $h$.
Observe that $w$ avoids the set
$\psi_{\varepsilon}(y)\cap\overline{B}(g_0,\alpha)$
because $f^{-1}(y)$ is of cotype $(q,\varepsilon,h_2(z))$ for all $z\in Q$.
This completes the proof of the claim.

\medskip
Let us go back to the proof of Proposition 4.8.
For a fixed $g\in C^*(X,l_2)$ and a function $\alpha\colon X\to (0,1)$
consider the set-valued map $\phi\colon Y\to C^*(X,l_2)$,
$\phi(y)=\overline{B}(g,\alpha)$, where $C^*(X,l_2)$
carries the uniform convergence topology.
According to the claim above, $\phi(y)\cap\psi_{\varepsilon}(y)$
is a $Z$-set in $\phi(y)$ for all $y\in Y$.
Moreover, $Y$ is a $C$-space, so we can apply \cite[Theorem 1.1]{gv}
to find a map $\theta\colon Y\to C^*(X,l_2)$ such that
$\theta(y)\in\overline{B}(g,\alpha)\backslash\psi_{\varepsilon}(y)$
for all $y\in Y$.
Finally, we define the map $\overline{g}\in C^*(X,l_2)$ by
$\overline{g}(x)=\theta(f(x))$, $x\in X$.
Then, $\overline{g}\in \overline{B}(g,\alpha)\cap\mathcal {F}_{\varepsilon}(y)$
for every $y\in Y$.
Hence, $\overline{g}\in \overline{B}(g,\alpha)\cap\mathcal {F}_{\varepsilon}$,
which completes the proof of Proposition 4.8.
\end{proof}

Therefore, the proof of Theorem 1.7 is done.
\section{Appendix}

Let us provide some more applications of our results.

\begin{thm}
Let $f\colon X\to Y$ be a perfect map between paracompact spaces
with $\dim Y=0$.
Suppose $\{F_i\}$ is a sequence of closed subsets of $X$ and $\{n_i\}$
a sequence of integers such that $\dim f|F_i\leq n_i$ for all $i$.
If $n\geq n_i$, $i=1,2,...$, then for every $m\geq n+1$ the space
$C^*(X,\R^m)$ contains a dense $G_{\delta}$-subset of maps $g$
such that for any $i$ and any $d$-plane $\Pi^d\subset\R^m$,
where $m-n_i\leq d\leq m$,
we have $\dim g\big(f^{-1}(y)\cap F_i\big)\cap\Pi^d\leq n_i+d-m$,
for any $y\in Y$.
\end{thm}

 The proof of Theorem 5.1 is based on Lemma 5.2 below.
Indeed, for every $i$ we apply Theorem 1.3 for the spaces $F_i$, $f(F_i)$
and the map $f|F_i$ to conclude that $C^*(F_i,\R^m)$ contains a dense
$G_{\delta}$-set $\mathcal{H}_i$ of maps $g$ such that
$\dim g\big(f^{-1}(y)\cap F_i\big)\cap\Pi^d\leq n_i+d-m$ for every
$y\in f(F_i)$ and every $d$-plane $\Pi^d\subset\R^m$ with $m-n_i\leq d\leq m$.
Since, by Lemma 5.2, each restriction map
$\pi_i\colon C^*(X,\R^m)\to C^*(F_i,\R^m)$ is open,
$\mathcal{K}_i=\pi_i^{-1}(\mathcal{H}_i)$ is dense and
$G_{\delta}$ in $C^*(X,\R^m)$.
Then the intersection of all $\mathcal{K}_i$
satisfies the requirements of Theorem 5.1.

\begin{lem}
Let $F$ be a closed subset of the normal space $X$ and $m\geq 1$.
Then the restriction map $\pi\colon C^*(X,\R^m)\to C^*(F,\R^m)$,
defined by $\pi(g)=g|F$, is open if both $C^*(X,\R^m)$ and $C^*(F,\R^m)$
are equipped simultaneously either with the source limitation topology
or the uniform convergence topology.
\end{lem}

Next corollary of Theorem 5.1 can be established in the same way
as Corollary 1.4 was obtained from Theorem 1.3.

\begin{cor}
Let $\{F_i\}$ be a sequence of closed subsets of the normal space $X$ with
$\dim F_i\leq n_i$.
If $m\geq n_i+1$ for each $i$, then $C^*(X,\R^m)$ equipped with
the uniform convergence topology contains a dense $G_{\delta}$-subset
of maps $g$ such that $\overline{g(F_i)}\cap\Pi^d$ is at
most $(n_i+d-m)$-dimensional for every $d$-plane $\Pi^d\subset\R^m$ with
$m-n_i\leq d\leq m$, $i=1,2,...$ .
\end{cor}

Our final application is an analogue of the Fox theorem \cite{f}
about economical extensions of maps.
If $A$ is a closed subset of a space $X$ and $h\in C^*(A,\R^m)$,
then $C^*_h(X,\R^m)$ denotes all maps $g\in C^*(X,\R^m)$ such that $g|A=h$.
Everywhere below $C^*_h(X,\R^m)$ is considered as a subspace of
$C^*(X,\R^m)$ with the uniform convergence topology.

\begin{cor}
Let $X$ be a normal space and $A$ a closed $G_{\delta}$-subset of  $X$
with $\dim (X\backslash A)\leq n$.
Then, for any $m\geq n+1$ and $h\in C^*(A,\R^m)$, there exists a dense and
$G_{\delta}$-set in $C^*_h(X,\R^m)$ of maps $g$ such that
$g(X\backslash A)\cap\Pi^d$ is at most $(n+d-m)$-dimensional
for every $d$-plane $\Pi^d\subset\R^m$ with $m-n\leq d\leq m$.
\end{cor}

\begin{proof}
Let $\{F_i\}$ be a sequence of closed subsets of $X$ with
$X\backslash A=\bigcup_{i=1}^{\infty}F_i$ and $h\in C^*(A,\R^m)$.
By Corollary 1.4, each $C^*(F_i,\R^m)$ contains a dense
$G_{\delta}$-subset $\mathcal{H}_i$ of maps $g$ with

\begin{itemize}
\item[(9)] $\dim \overline{g(F_i)}\cap\Pi^d\leq n+d-m$ for any
$d$-plane $\Pi^d\subset\R^m$, where $m-n\leq d\leq m$.
\end{itemize}

\noindent
Using that $F_i$ and $A$ are closed and disjoint subsets of $X$,
one can show that every projection
$p_i\colon C^*_h(X,\R^m)\to C^*(F_i,\R^m)$, $p_i(g)=g|F_i$,
is open and surjective when both $C^*_h(X,\R^m)$ and $C^*(F_i,\R^m)$
are equipped with the uniform convergence topology.
Hence, $\mathcal{K}_i=p_i^{-1}(\mathcal{H}_i)$ is dense and
$G_{\delta}$ in $C^*_h(X,\R^m)$.
Since $C^*_h(X,\R^m)$ has the Baire property,
$\mathcal{K}=\bigcap_{i=1}^{\infty}\mathcal{K}_i$ is dense and
$G_{\delta}$ in $C^*_h(X,\R^m)$.
Finally, observe that (9) implies $\dim g(X\backslash A)\cap\Pi^d\leq n+d-m$
for every $g\in\mathcal{K}$ and every $d$-plane $\Pi^d\subset\R^m$
with $m-n\leq d\leq m$.
\end{proof}

Finally, let us discuss some possible improvements of our results.

\smallskip\noindent
{\bf Conjecture 1.}
{\it Let $f\colon X\to Y$ be a map of finite dimensional metrizable compacta.
Then, $C(X,\R^m)$ contains a dense $G_{\delta}$-subset of maps
$\varphi$ such that, for any integers $d,t,T$ with $0\leq t\leq d\leq T\leq m$
and $\dim f+d+1\leq m$ and for any $d$-plane
$\Pi^d\subset\R^{m}$, parallel to some coordinate planes $\Pi^t\subset\Pi^T$ in $\R^m$,
each set
$f^{-1}(y)\cap \varphi^{-1}(\Pi^d)$, $y\in Y$, has no more than
$$
1+\frac{\dim Y+\dim f+(T-d)(d-t)}{m-\dim f-d}
$$
points.}

\medskip
Obviously, Theorem 1.5 implies the validity of Conjecture 1 in case
$Y$ is $0$-dimensional.

\medskip\noindent
{\bf Conjecture 2.}
{\it Let $f\colon X\to Y$ be a map of finite-dimensional compacta.
Then, $C(X,\R^m)$ contains a dense $G_{\delta}$-subset
of maps $\varphi$ such that
$$
\dim \bigl(\varphi \bigl(f^{-1}(y)\bigr)\cap \Pi^d\bigr)\leq \dim f+d-m,
$$
for any $d$-plane $\Pi^d\subset \R^{m}$ with $m-\dim f\leq d\leq m$
and any $y\in Y$.}

\medskip
Theorem 1.3 provides the validity of Conjecture 2 in case
$Y$ is $0$-dimensional.
As in Theorem 1.3, it suffices to establish Conjecture 2
in the special case $d=m-\dim f$.
Conjecture 2 is also true when $\dim f=0$, this is Uspenskij's theorem \cite{u}
about light mappings (see also \cite{tv1}
for a generalization of the Uspenskij result when $\dim f>0$).

Let $f\colon X\to Y$, $\varphi\in C(X,\R^m)$ and $t,d,T$ be integers with
$0\leq t\leq d\leq T\leq m$ and $d-t+1\leq q$.
Below we consider the set $\displaystyle B_{q,d,t,T}^f(\varphi)$
consisting of all points $\displaystyle (y,y_1,...,y_q)\in Y\times (\R^m)^q$
satisfying the following condition:
there exist points $x_1,...,x_q\in f^{-1}(y)$ with
$x_i\ne x_j \text{ for } i\ne j$ such that $y_i=\varphi (x_i)$
and all $y_i$, $i=1,...,q$, belong to a $d$-plane in $\R^m$
parallel to some coordinates planes $\Pi^t\subset\Pi^T\subset\R^m$.

\medskip\noindent
{\bf Conjecture 3.}
{\it Let $f\colon X\to Y$ be a map of finite-dimensional metrizable compacta.
Then $C(X,\R^m)$ contains a dense $G_{\delta}$-set $\mathcal H$ of maps
$\varphi$ such that
$$
\dim B_{q,d,t,T}^f(\varphi )\leq \dim Y+\dim f+(T-d)(d-t)-(q-1)(m-\dim f-d)
$$
for any integers $d,t,T,q$ satisfying the conditions $0\leq t\leq d\leq T\leq m$,
$d-t+1\leq q$ and $\dim f+d+1\leq m$.}

\medskip
If the right side of the inequality from Conjecture 3 is $\leq -1$,
then the set $\displaystyle B_{q,d,t,T}^f(\varphi )$ is empty.
Conditions when $\displaystyle B_{q,d,t,T}^f(\varphi )$ is empty
are discussed in Conjecture 1.
If $d=0$, the set $\displaystyle B_{q,0,0,T}^f(\varphi )$
does not depend on $T$.
In this case it is homeomorphic to the set
$\displaystyle B_{q}^f(\varphi)=\{(y,z)\in Y\times\R^m:
|f^{-1}(y)\cap\varphi^{-1}(z)|\geq q\}$.
For $d=0$ and $Y$ a point, the statement of Conjecture 3 was obtained
in 1933 by Hurewicz \cite{h}
(recall that the maps $\varphi\in\mathcal H$ from the Hurewicz theorem
are called regularly branched maps \cite{drs}).
A parametric version of the Hurewicz result
was obtained in \cite[Theorem 1.1]{tv2}.
Following the terminology of \cite{tv2}, the inequality from Conjecture 3
is satisfied for $d=0$ and all $q\geq 1$ if and only if $\varphi$
is an $f$-regularly branched map.
By \cite[Theorem 1.1]{tv2}, $C(X,\R^m)$ contains a dense $G_{\delta}$-subset
of $f$-regularly branched maps, so Conjecture 3 holds for $d=0$.
Let us note that for $d\geq 1$, even in the case when $Y$ is a point,
the conjecture is open.

For maps $\displaystyle f_i\colon X_i\to \R^m$, $i=1,...,q$, and
integers $0\leq t\leq d\leq T\leq m$, let
$\displaystyle B_{d,t,T}(f_1,...,f_q)$ be the set of the points
$\displaystyle (y_1,...,y_q)\in(\R^m)^q$ such that the points
$\displaystyle y_i=f(x_i)$, $i=1,...,q$, belong to a $d$-plane in $\R^m$
parallel to some coordinate planes $\Pi^t\subset\Pi^T$.
Let also
$$
C_{d,t,T}(f_1,...,f_q)=
\bigl\{(x_1,..,x_q)\in X_1\times ..\times X_q:
\bigl(f(x_1),..,f(x_q)\bigr)\in B_{d,t,T}(f_1,..,f_q)\bigr\}.
$$

\medskip\noindent
{\bf Conjecture 4.}
{\it Let the numbers $n_1,...,n_q,m,d,t,T$ satisfy the inequalities
$0\leq t\leq d\leq T\leq m$, $0\leq n_1$, ..., $0\leq n_q$,
$n_1+1+d\leq m$, ..., $n_q+1+d\leq m$, $d-t+1\leq q$ and
$$
n_1+...+n_q\geq (m-d)(q-1)-(T-d)(d-t).
$$
Then, there exist $\varepsilon >0$ and maps
$\displaystyle f_i\colon\Delta^{n_i}\to \R^m$, $i=1,...,q$,
such that the set $\displaystyle B_{d,t,T}(g_1,...,g_q)$ is not empty
for any maps $\displaystyle g_i\colon\Delta^{n_i}\to \R^m$
with each $g_i$ being $\varepsilon$-close to $f_i$, $i=1,...,q$.}

In the simplest case ($q=1$, $n_1=0$, $n_2=m$, $d=0$),
the statement of Conjecture 4 is exactly the Alexandrov
theorem that the identity map of the ball onto itself is essential.
For $d=0$ Conjecture 4 is also true
(see for example \cite[Corallary 3]{b1}, \cite[Lemma 6.3, p.65]{bms}).
In the case $d=q-1$, $t=0$ nd $T=m$ Conjecture 4 was proved
by Boltyanski \cite[Lemma 9]{b} (this is the main ingredient in his proof that
$m\geq nk+n+k$ provided for some $n$-dimensional polyhedron $X$ the set of
$k$-regular maps from $X$ into $\R^m$ is dense in $C(X,\R^m)$).

\medskip\noindent
{\bf Conjecture 5.}
{\it Let the numbers $n_1,...,n_q,m,d,t,T$ satisfy the inequalities
$0\leq t\leq d\leq T\leq m$, $0\leq n_1$, ..., $0\leq n_q$,
$n_1+1+d\leq m$, ..., $n_q+1+d\leq m$, $d-t+1\leq q$.
Then, there exist $\varepsilon >0$ and maps
$f_i:\Delta^{n_i}\to \Bbb R^m$, $i=1,...,q$, such that
$$
\dim C_{d,t,T}(g_1,...,g_q)\geq n_1+...+n_q-(m-d)(q-1)+(T-d)(d-t)
$$
for any maps $g_i:\Delta^{n_i}\to \Bbb R^m$ with each
$g_i$ being $\varepsilon$-close to $f_i$, $i=1,...,q$}.

For $d=0$ the statement of Conjecture 5
was established in \cite[Corallary 3]{b1}.

We are going to finish with the following problem,
where $\displaystyle \Delta_n^N$ denotes the $n$-dimensional skeleton of
the $N$-dimensional simplex $\displaystyle\Delta^N$.

\medskip\noindent
{\bf Problem.}
{\it Find all integers $n,m,q,d,t,T$
with $0\leq t\leq d\leq T\leq m$, $n+1+d\leq m$, $d-t+1\leq q$,
$$
n\geq (m-n-d)(q-1)-(T-d)(d-t)
\eqno (n,m,q,d,t,T)
$$
and satisfying the following condition:
There exists a number $N$ such that for any mapping
$\displaystyle f\colon\Delta_n^N\to \R^m$ we can find pairwise disjoint
simplexes $\displaystyle\sigma_1,...,\sigma_q\subset \Delta_n^N$
whose images under $f$ meet a $d$-plane $\Pi^d\subset\R^m$
parallel to some coordinate planes $\Pi^t\subset\Pi^T$}.

Observe that, if the above problem has a positive solution for some numbers
$n,m,q,d,t,T$ satisfying the assumptions of the problem, then Theorem 2
can not be improved, not only on the level of a dense set of mappings,
but also on the level of {\it existence} of one map with preimages having small
cardinality, even in the class of polyhedra.
Most of the results providing the existence of such a number $N$
were established for $d=0$:
van Kampen and Flores \cite{k}, \cite{fl} ($q=2$ and $N=2n+2$);
Sarkaria \cite{s} ($q$ prime and $N=qn+2q-2$);
Volovikov \cite{v} ($q$ being a power of a prime number and $N=qn+2q-2$);
Bogatyi \cite[Corollary 11]{b1} ($q=n+1$ and $N\leq 2n^2+5n$).
Another result of Bogatyi \cite{b4} provides such $N$
for $d=q-1$, $t=0$, $T=m$ with $q$ being odd.
\v{Z}ivaljevi\'c obtained \cite{z} a result about emdeddings of the graph
$K_{6,6}$ into $\R^3$ which implies a positive solution of the problem for
$n=1$, $m=3$, $t=0$, $T=3$, $q=4$ and $N=11$.

If $d=q-1$, $q=2$, $t=0$ and $T=m$, there is no any number $N$ satisfying the
conditions in the above problem, see \cite{b}.

The above problem can have a little different treatment:
{\it For a given family of integers $n,q,d,t,T$ find the biggest number $m$
for which such number $N$ does exist}.
Then the question to determine the smallest such number $N$ is arising. In 
this case, more complicated is the description of such minimal 
subpolyhedra in $\Delta_n^N$.
For $n=1$, $d=0$ and $m=2$ the Kuratowski graphs 
$K_5$ and $K_{3,3}$ are minimal subpolyhedra.

Let us finally note that there exist connections between the above problem and
conjectures about different forms of the Tverberg theorem \cite{z}, \cite{b3}.
It is also well known that $k$-regular mappings are closely related to
interpolation and approximation problems.
In view of this, it is important to find applications of the mappings described
by Theorem 1.5 and Theorem 1.7 in interpolation and approximation problems.


\bigskip
\end{document}